\documentclass{amsart}

\newtheorem{theorem}{Theorem}[section]
\newtheorem{lemma}[theorem]{Lemma}
\newtheorem{proposition}[theorem]{Proposition}

\theoremstyle{definition}
\newtheorem{definition}[theorem]{Definition}

\newtheorem{corollary}[theorem]{Corollary}
\theoremstyle{remark}
\newtheorem{remark}[theorem]{Remark}

\numberwithin{equation}{section}

\begin{document}

\title{Cappell-Miller analytic torsion for manifolds with boundary}

\author{Rung-Tzung Huang}
\address{Institute of Mathematics, Academia Sinica, Nankang 11529, Taipei, Taiwan}

\email{rthuang@math.sinica.edu.tw}



\begin{abstract}
Inspired by the work of Boris Vertman on refined analytic torsion for manifolds with boundary, in this paper we extend the construction of the Cappell-Miller analytic torsion to manifolds with boundary. We also compare it with the refined analytic torsion on manifolds with boundary. As a byproduct of the gluing formula for refined analytic torsion and the comparison theorem for the Cappell-Miller analytic torsion and the refined analytic torsion, we establish the gluing formula for the Cappell-Miller analytic torsion in the case that the Hermitian metric is flat. \end{abstract}

\maketitle

\section{Introduction}
Let $E$ be a flat complex vector bundle over a closed oriented odd dimensional manifold $M$. Braverman and Kappeler \cite{BK1,BK3,BK4,BK5,BK6} defined and studied the refined analytic torsion for $(M,E)$, which can be viewed as a refinement of the Ray-Singer torsion \cite{RS} and an analytic analogue of the Farber-Turaev torsion, \cite{FT1,FT2,Tu1,Tu2}. It was shown that the refined analytic torsion is closely related with the Farber-Turaev torsion, \cite{BK1,BK3,BK6,H}. Burghelea and Haller \cite{BH1,BH2,BH3} defined the complex valued Ray-Singer torsion associated to a non-degenerate symmetric bilinear form on a flat vector bundle over an arbitrary dimensional manifold and make an explicit conjecture between the Burghelea-Haller analytic torsion and the Farber-Turaev torsion. This conjecture was proved up to sign by Burghelea-Haller \cite{BH3} and in full generality by Su-Zhang \cite{SZ}. Cappell and Miller \cite{CM} used non-self-adjoint Laplace operators to define another complex valued analytic torsion and used the method in \cite{Su} to prove an extension of the Cheeger-M\"{u}ller theorem, \cite{Ch,Mu,Mu2}, which states that the Cappell-Miller analytic torsion is equal to the Reidemeister torsion of the bundle $E\oplus E^*$, where $E^*$ denotes the dual bundle to $E$.
 
By combining the absolute and relative boundary conditions, Vertman \cite{Ve1} applied the original construction of Braverman-Kappeler \cite{BK1,BK3} to a new setting. The proposed construction refines the square of the Ray-Singer torsion, but applies to compact manifolds with and without boundary. In a subsequent paper \cite{Ve2} Vertman derived a gluing formula for the refined analytic torsion in this setting under the assumption that the Hermitian metric is flat. Inspired by the construction of \cite{Ve1}, Su \cite{Su} extended the Burghelea-Haller analytic torsion to compact manifolds with boundary and compared it with the refined analytic torsion. In this paper we extend the construction of the Cappell-Miller analytic torsion to manifolds with boundary and compare the Cappell-Miller analytic torsion with the refined analytic torsion. As a byproduct of the gluing formula for refined analytic torsion \cite{Ve2} and the comparison theorem for the Cappell-Miller analytic torsion and the refined analytic torsion, we establish the gluing formula for the Cappell-Miller analytic torsion in the case that the Hermitian metric is flat. It would be interesting to establish the gluing formula for the Cappell-Miller analytic torsion in the case that the Hermitian metric is not necessarily flat.

The rest of this paper is organized as follows. In Section 2, we extend the construction of the Cappell-Miller analytic torsion to manifolds with boundary. In Section 3, we compare the extended Cappell-Miller analytic torsion with the refined analytic torsion. In Section 4, we establish the gluing formula for the Cappell-Miller analytic torsion in the case that the Hermitian metric is flat.

\subsection*{Acknowledgement}
The author would like to thank Maxim Braverman for helpful comments.
 
\section{The Cappell-Miller analytic torsion for manifolds with boundary}

Inspired by the paper \cite{Ve1}, in this section we generalize the construction of the Cappell-Miller analytic torsion to manifolds with boundary.

\subsection{The Cappell-Miller torsion for finite dimensional complexes} Given a complex vector space $V$ of dimension $l$, the {\em determinant line} of $V$ is the line $\operatorname{Det}(V):=\wedge^l V$, where $\wedge^l V$ denotes the $l$-th exterior power of $V$. By definition, we set $\operatorname{Det}(0):=\mathbb{C}$. Further, we denote by $\operatorname{Det}(V)^{-1}$ the dual line of $\operatorname{Det}(V)$. Let 
\begin{equation}\label{E:detline} 
 (C^\bullet,d) \  : \ 0 \stackrel{d}{\longrightarrow} \ C^0 \ \stackrel{d}{\longrightarrow}\ C^1\ \stackrel{d}{\longrightarrow} \ \cdots \stackrel{d}{\longrightarrow}\ \ C^n \  \stackrel{d}{\longrightarrow} 0
\end{equation}
be a cochain complex of finite dimensional complex vector spaces. Denote by $H^\bullet(d)=\bigoplus_{i=0}^n H^i(d)$ its cohomology. Set
\begin{equation}
\operatorname{Det}(C^\bullet)\,:=\,\bigotimes_{j=0}^n \operatorname{Det}(C^j)^{(-1)^j} , \qquad \operatorname{Det}(H^\bullet(d))\,:=\, \bigotimes_{j=0}^n \operatorname{Det}(H^j(d))^{(-1)^j}.
\end{equation}
There is a standard isomorphism, cf. \cite[(6.9)]{CM},
\begin{equation}\label{E:isomorphism}
\tau \,:\,\operatorname{Det}(C^\bullet) \rightarrow \operatorname{Det}(H^\bullet(d)).
\end{equation}

Now if in addition, $C^\bullet$ has another differential $d^* :C^{j} \to C^{j-1}, (d^*)^2=0$, then
\[
(C^\bullet,d^*) \  : \ 0 \stackrel{d^*}{\longleftarrow} \ C^0 \ \stackrel{d^*}{\longleftarrow}\ C^1\ \stackrel{d^*}{\longleftarrow} \ \cdots \stackrel{d^*}{\longleftarrow}\ \ C^n \  \stackrel{d^*}{\longleftarrow} 0
\]
is a chain complex of finite dimensional complex vector spaces.
Denote by $H_\bullet(d^*)=\bigoplus_{i=0}^n H_i(d^*)$ its homology.
Set
\[
\operatorname{Det}(H_\bullet(d^*))\,:=\, \bigotimes_{j=0}^n \operatorname{Det}(H_j(d^*))^{(-1)^j}.
\]
Then correspondingly there is a standard isomorphism, cf. \cite[(6.16)]{CM},
\begin{equation}\label{E:isomorphism2}
\tau' \,:\,\operatorname{Det}(C^\bullet) \rightarrow \operatorname{Det}(H_\bullet(d^*)).
\end{equation}

Let $c_{j}\in \operatorname{Det}(C^{j})$ $(j=0,\cdots,n)$ and denote by $c_{j}^{-1}$ the unique element of $\operatorname{Det}(C^{j})^{-1}$ such that $c^{-1}_{j}(c_{j})=1$. Consider the element
\[
c\, := \, c_{0} \otimes c_{1}^{-1} \otimes \cdots \otimes c_{n}^{(-1)^n} 
\]
of $\operatorname{Det}(C^\bullet)$. Then for the bi-graded complex $(C^\bullet, d, d^*)$, the Cappell-Miller torsion is the algebraic torsion invariant
\begin{equation}\label{E:tautor}
\tau(C^\bullet, d, d^*) \, := \, (-1)^{S(C^\bullet)}\tau(c) \otimes (\tau'(c))^{-1} \,  \in \, (\operatorname{Det}(H^\bullet(d)) )\otimes (\operatorname{Det}(H_\bullet(d^*)))^{-1},
\end{equation} 
where $(-1)^{S(C^\bullet)}$ is defined in \cite[Section 6]{CM}. 

\begin{remark}
Note that Braverman and Kappeler \cite[(2-14)]{BK3} introduce a sign refined version of the standard isomorphism $\tau$ or $\tau'$, cf. \cite{CM,Mi}, to obtain various compatibility properties. Instead of this, Cappell and Miller introduce the total sign correction $(-1)^{S(C^\bullet)}$ in the above definition of $\tau(C^\bullet,d,d^*)$ to reestablish desirable compatibility properties.     
\end{remark}

The following proposition was proved in \cite{CM}.

\begin{proposition}
Suppose that $(C^\bullet,d,d^*)$ is a finite bi-graded complex and the combinatorial Laplacian, $\Delta_{j} \, := \, (d+d^*)^2|_{C^j}$ $(j=0,\cdots,n)$, has no zero eigenvalue. Then the cohomology groups $H^{j}(d)$ and the homology groups $H_j(d^*)$ vanish and 
\[
\tau(C^\bullet,d,d^*) \, = \, \prod_{j=0}^n(\det(\Delta_j))^{{(-1)}^{j+1}j}.
\]   
\end{proposition}


\subsection{Fredholm complexes for compact manifolds}
Let $(M,g^M)$ be a compact oriented Riemannian manifold with boundary $\partial M$, possibly empty, where $g^M$ is the Riemannian metric on $M$. Suppose that $E$ is a complex vector bundlle over $M$ endowed with a flat connection $\nabla$. The connection $\nabla$ gives rise to a covariant differential on $\Omega^\bullet_0(M,E)$, the space of smooth $E$-valued differential forms with compact support in the interior of the manifold $M$. We choose a Hermitian metric $h^E$ so that together with the Riemannian metric $g^M$ we can define an $L^2$-inner product $<\, , \,>_M$ on $\Omega^\bullet_0(M,E)$. Denote by $L^2_\bullet(M,E)$ the $L^2$-completion of $\Omega^\bullet_0(M,E)$. We choose the dual connection $\nabla^{\prime}$ with respect to
$h^{E}$ satisfying the following property. For $\phi$, $\psi \in
C^{\infty}(M, E)$,
$$
d ( h^{E}(\phi, \psi) ) = h^{E}(\nabla\phi, \psi) + h^{E}(\phi,
\nabla^{\prime}\psi).
$$
We then extend $\nabla^{\prime}$ to a covariant differential
$\nabla^{\prime} : \Omega^{\bullet}_0(M, E) \rightarrow \Omega^{\bullet + 1}_0(M, E)$.
Consider the differential operators $\nabla$, $\nabla'$ and their formal adjoint differential operators $\nabla^t$, $\nabla'^t$. The associated minimal closed extensions $\nabla_{\operatorname{min}}$, $\nabla_{\operatorname{min}}'$ and $\nabla_{\operatorname{min}}^t$, $\nabla_{\operatorname{min}}'^t$ are defined as the graph-closures in $L^2_\bullet(M,E)$ of the corresponding differential operators. The maximal closed extensions are defined by 
\[
\nabla_{\operatorname{max}}:=(\nabla_{\operatorname{min}}^t)^*, \qquad \nabla_{\operatorname{max}}':=(\nabla_{\operatorname{min}}'^t)^*.
\]
These extensions define Hilbert complexes in the following sense, as introduced in \cite{BL}.

\begin{definition}[\cite{BL}]
Let the Hilbert spaces $H_i$, $i=0,\cdots,m$, $H_{m+1}= \{ 0 \}$ be mutually orthogonal. For each $i=0,\cdots,m$, let $D_i \in C(H_i,H_{i+1})$ be a closed operator with domain $\mathcal{D}(D_i)$ dense in $H_i$ and range in $H_{i+1}$. Put $\mathcal{D}_i=\mathcal{D}(D_i)$ and $R_i=D_i(\mathcal{D}_i)$, and assume
\[
R_i \subseteq \mathcal{D}_{i+1}, \quad D_{i+1} \circ D_i = 0. 
\]
This defines a complex $(\mathcal{D},D)$
\[
  0 \longrightarrow \ \mathcal{D}_0 \ \stackrel{D_0}{\longrightarrow}\ \mathcal{D}_1 \ \stackrel{D_1}{\longrightarrow} \ \cdots \stackrel{D_{m-1}}{\longrightarrow}\ \ \mathcal{D}_m \ \longrightarrow 0.
\]
Such a complex is called a Hilbert complex. If the homology of the complex is finite, i.e. if $R_i$ is closed and $\ker D_i/ \operatorname{im} D_{i-1}$ is finite-dimensional for all $i=0,\cdots,m$, the complex is referred to as a Fredholm complex. 
\end{definition}

By \cite[Lemma 3.1]{BL} we have Hilbert complexes $(\mathcal{D}_{\operatorname{min}},\nabla_{\operatorname{min}})$ and $(\mathcal{D}_{\operatorname{max}},\nabla_{\operatorname{max}})$, where $\mathcal{D}_{\operatorname{min}}=\mathcal{D}(\nabla_{\operatorname{min}})$ and $\mathcal{D}_{\operatorname{max}}=\mathcal{D}(\nabla_{\operatorname{max}})$.
The Laplace operators, associated to the Hilbert complexes $(\mathcal{D}_{\operatorname{min}},\nabla_{\operatorname{min}})$ and $(\mathcal{D}_{\operatorname{max}},\nabla_{\operatorname{max}})$, are respectively defined as follows:
\[
\Delta_{\operatorname{rel}} \, = \, (\nabla_{\operatorname{min}} + \nabla_{\operatorname{min}}^*)^2,
\]
\[
\mathcal{D}(\Delta_{\operatorname{rel}}) \, = \, \{ \omega \in \mathcal{D}(\nabla_{\operatorname{min}}) \cap \mathcal{D}(\nabla_{\operatorname{min}}^*) \, | \, \nabla_{\operatorname{min}} \omega \in \mathcal{D}(\nabla_{\operatorname{min}}^*), \, \nabla_{\operatorname{min}}^* \omega \in \mathcal{D}(\nabla_{\operatorname{min}}) \}
\]
and
\[
\Delta_{\operatorname{abs}} \, = \, (\nabla_{\operatorname{max}} + \nabla_{\operatorname{max}}^*)^2,
\]
\[
\mathcal{D}(\Delta_{\operatorname{abs}}) \, = \, \{ \omega \in \mathcal{D}(\nabla_{\operatorname{max}}) \cap \mathcal{D}(\nabla_{\operatorname{max}}^*) \, | \, \nabla_{\operatorname{max}} \omega \in \mathcal{D}(\nabla_{\operatorname{max}}^*), \, \nabla_{\operatorname{max}}^* \omega \in \mathcal{D}(\nabla_{\operatorname{max}}) \}.
\]
The following theorem \cite[Theorem 3.2]{Ve1} is the twisted setup of \cite[Theorem 4.1]{BL}.
\begin{theorem}\label{T:hilbert}
The Hilbert complexes $(\mathcal{D}_{\operatorname{min}},\nabla_{\operatorname{min}})$ and $(\mathcal{D}_{\operatorname{max}},\nabla_{\operatorname{max}})$ are Fredholm with the associated Laplacians $\Delta_{\operatorname{rel}}$ and $\Delta_{\operatorname{abs}}$ being strongly elliptic in the sense of \cite[Subsection 1.11]{Gi}. The de Rham isomorphism identifies the homology of the complexes with the relative and absolute cohomology with coefficients
$$H^\bullet_{\operatorname{rel}}(M,E) := H^\bullet(M,\partial M,E) \cong H^\bullet(\mathcal{D}_{\operatorname{min}},\nabla_{\operatorname{min}}), 
$$
$$
H^\bullet_{\operatorname{abs}}(M,E) := H^\bullet(M,E) \cong H^\bullet(\mathcal{D}_{\operatorname{max}},\nabla_{\operatorname{max}}).
$$
Furthermore the cohomology of Fredholm complexes $(\mathcal{D}_{\operatorname{min}},\nabla_{\operatorname{min}})$ and $(\mathcal{D}_{\operatorname{max}},\nabla_{\operatorname{max}})$ can be computed from the following smooth subcomplexes
$$
(\Omega^\bullet_{\operatorname{min}}(M,E),\nabla), \  \ \Omega^\bullet_{\operatorname{min}}(M,E) \, = \, \{ \, \omega \in \Omega^\bullet(M,E) | \iota^*(\omega)=0 \, \},$$
$$
(\Omega^\bullet_{\operatorname{max}}(M,E),\nabla), \ \ \Omega^\bullet_{\operatorname{max}}(M,E) \, = \, \Omega^\bullet(M,E),
$$
respectively, where we denote by $\iota: \partial M \to M$ the natural inclusion of the boundary.
\end{theorem}

\subsection{Non-self-adjoint Laplacian operators}
The Riemannian metric $g^M$ and the fixed orientation on $M$ give rise to the Hodge star operator $\star$, which induces an isomorphism on the spaces of forms $\Omega^\bullet(M,E)$ and extends to $L^2_\bullet(M,E)$, also denote by $\star$. Using the Hodge star operator $\star$, we define the involution $\Gamma=\Gamma(g^M):\Omega^\bullet(M,E) \to \Omega^{n-\bullet}(M,E)$ by
\[
\Gamma \omega \, := \, i^r(-1)^{\frac{p(p+1)}{2}}\star \omega, \quad \omega \in \Omega^p(M,E),
\]
where $r=\frac{n+1}{2}$ if $n$ is odd and $r=\frac{n}{2}$ if $n$ is even. It is straightforward to see that $\Gamma^2=\operatorname{Id}$. We recall the following lemma, cf. \cite[Lemma 3.3]{Ve1}.
\begin{lemma}
The chirality operator $\Gamma$ on $L_\bullet^2(M,E)$, restricted to $\mathcal{D}(\nabla_{\operatorname{min}})$ and $\mathcal{D}(\nabla_{\operatorname{max}})$, acts as follows:
\[
\Gamma |_{\mathcal{D}(\nabla_{\operatorname{min}})} \, : \, \mathcal{D}(\nabla_{\operatorname{min}}) \to \mathcal{D}((\nabla_{\operatorname{max}}')^*),
\]
\[
\Gamma |_{\mathcal{D}(\nabla_{\operatorname{max}})} \, : \, \mathcal{D}(\nabla_{\operatorname{max}}) \to \mathcal{D}((\nabla_{\operatorname{min}}')^*).
\]
When $\Gamma$ restricted to appropriate domains, we have
\begin{equation}\label{E:chirality}
\begin{array}{l}
\nabla_{\operatorname{max}}^\sharp  :=  (\nabla_{\operatorname{max}}')^*  =  \Gamma|_{\mathcal{D}(\nabla_{\operatorname{min}})} \nabla_{\operatorname{min}} \Gamma|_{\mathcal{D}(\nabla_{\operatorname{max}}^\sharp)},  \\ \nabla_{\operatorname{min}}^\sharp  :=  (\nabla_{\operatorname{min}}')^*  =  \Gamma|_{\mathcal{D}(\nabla_{\operatorname{max}})} \nabla_{\operatorname{max}} \Gamma|_{\mathcal{D}(\nabla_{\operatorname{min}}^\sharp)}.
\end{array}
\end{equation}
\end{lemma} 
It is easy to check that $(\nabla_{\operatorname{max}}^\sharp)^2=0$ and $(\nabla_{\operatorname{min}}^\sharp)^2=0$. The Hermitian metric $h^E$ defines a conjugate linear bundle isomorphism from $E$ to its dual $E^*$, also denote by $h^E$. Then the adjoints of $\nabla_{\operatorname{max}}$ and $\nabla_{\operatorname{min}}$ are 
\[
\nabla_{\operatorname{max}}^* \, = \, \Gamma|_{\mathcal{D}(\nabla_{\operatorname{min}})} \, (h^E)^{-1} \, \nabla_{\operatorname{min}} \, h^E \, \Gamma|_{\mathcal{D}(\nabla_{\operatorname{max}}^*)} 
\]
and
\[
\nabla_{\operatorname{min}}^* \, = \,\Gamma|_{\mathcal{D}(\nabla_{\operatorname{max}})} \,  (h^E)^{-1} \, \nabla_{\operatorname{max}} \, h^E \, \Gamma|_{\mathcal{D}(\nabla_{\operatorname{min}}^*)},
\]
respectively. The operator $\nabla_{\operatorname{min}}^*$ (resp. $\nabla_{\operatorname{max}}^*$) differs from the operator $\nabla_{\operatorname{min}}^\sharp$ (resp. $\nabla_{\operatorname{max}}^\sharp$) by an endomorphism-valued differential form of degree one, which can be viewed as a differential operator of degree zero. A differential operator of degree zero naturally extends to a bounded operator on the $L^2$-Hilbert space, and hence does not pose additional restrictions on the domain, \cite[P. 1995]{Ve1}. Therefore,
\begin{equation}\label{E:doeq}
\mathcal{D}(\nabla_{\operatorname{min}}^*) = \mathcal{D}(\nabla_{\operatorname{min}}^\sharp),  \quad \mathcal{D}(\nabla_{\operatorname{max}}^*) = \mathcal{D}(\nabla_{\operatorname{max}}^\sharp). 
\end{equation}
We define the flat Laplace operators, associated to Hilbert complexes $(\mathcal{D}_{\operatorname{min}},\nabla_{\operatorname{min}})$ and $(\mathcal{D}_{\operatorname{max}},\nabla_{\operatorname{max}})$, respectively as follows:  
\[
\Delta_{\operatorname{rel}}^\sharp \, = \, (\nabla_{\operatorname{min}} + \nabla_{\operatorname{min}}^\sharp)^2,
\]
\[
\mathcal{D}(\Delta_{\operatorname{rel}}^\sharp) \, = \, \{ \omega \in \mathcal{D}(\nabla_{\operatorname{min}}) \cap \mathcal{D}(\nabla_{\operatorname{min}}^\sharp) \, | \, \nabla_{\operatorname{min}} \omega \in \mathcal{D}(\nabla_{\operatorname{min}}^\sharp), \, \nabla_{\operatorname{min}}^\sharp \omega \in \mathcal{D}(\nabla_{\operatorname{min}}) \}
\]
and
\[
\Delta_{\operatorname{abs}}^\sharp \, = \, (\nabla_{\operatorname{min}} + \nabla_{\operatorname{max}}^\sharp)^2,
\]
\[
\mathcal{D}(\Delta_{\operatorname{abs}}^\sharp) \, = \, \{ \omega \in \mathcal{D}(\nabla_{\operatorname{max}}) \cap \mathcal{D}(\nabla_{\operatorname{max}}^\sharp) \, | \, \nabla_{\operatorname{max}} \omega \in \mathcal{D}(\nabla_{\operatorname{max}}^\sharp), \, \nabla_{\operatorname{max}}^\sharp \omega \in \mathcal{D}(\nabla_{\operatorname{max}}) \}.
\]
They are generally non-self-adjoint operators, as opposed to the standard self-adjoint Laplace operators $\Delta_{\operatorname{rel}}$ and $\Delta_{\operatorname{abs}}$. However, $\Delta_{\operatorname{rel}}$ (resp. $\Delta_{\operatorname{abs}}$) and $\Delta_{\operatorname{rel}}^\sharp$ (resp. $\Delta_{\operatorname{abs}}^\sharp$) have the same leading symbols. Hence the spectrum of the operators $\Delta_{\operatorname{rel}}^\sharp$ and $\Delta_{\operatorname{abs}}^\sharp$ are discrete, and, by \eqref{E:doeq},
\[
\mathcal{D}(\Delta_{\operatorname{rel}}^\sharp) \, = \, \mathcal{D}(\Delta_{\operatorname{rel}}), \quad \mathcal{D}(\Delta_{\operatorname{abs}}^\sharp) \, = \, \mathcal{D}(\Delta_{\operatorname{abs}}).
\] 
Moreover, if the Hermitian metric $h^E$ is flat, then $(h^E)^{-1}\nabla_{\operatorname{min}} h^E = \nabla_{\operatorname{min}}$ (resp. $(h^E)^{-1}\nabla_{\operatorname{max}} h^E = \nabla_{\operatorname{max}}$) and, hence, $\Delta_{\operatorname{rel}}^\sharp = \Delta_{\operatorname{rel}}$ (resp. $\Delta_{\operatorname{abs}}^\sharp = \Delta_{\operatorname{abs}}$). 

\subsection{The Cappell-Miller analytic torsion for manifolds with boundary}
We now discuss the construction of the Cappell-Miller analytic torsion for the complex $(\mathcal{D}_{\operatorname{min}},\nabla_{\operatorname{min}},\nabla_{\operatorname{min}}^\sharp)$ explicitly. The construction of the Cappell-Miller analytic torsion for the complex $(\mathcal{D}_{\operatorname{max}},\nabla_{\operatorname{max}},\nabla_{\operatorname{max}}^\sharp)$ is exactly the same.

Let $\lambda \ge 0$ be any nonnegative real number. Denote by $\Pi_{\Delta_{\operatorname{rel}}^\sharp,[0,\lambda]}$ the spectral projection of $\Delta_{\operatorname{rel}}^\sharp$ onto the eigenspaces of absolute value in $[0,\lambda]$:
\begin{equation}\label{E:specpro}
\Pi_{\Delta_{\operatorname{rel}}^\sharp,[0,\lambda]} := \frac{i}{2\pi}\int_{C(\lambda)}(\Delta_{\operatorname{rel}}^\sharp-x)^{-1}dx,
\end{equation} 
where $C(\lambda)$ being any closed counterclockwise circle surrounding eigenvalues of absolute value in $[0,\lambda]$ with no other eigenvalue inside. The dimension of the image of $\Pi_{\Delta_{\operatorname{rel}}^\sharp,[0,\lambda]}$ is finite. In particular, $\Pi_{\Delta_{\operatorname{rel}}^\sharp,[0,\lambda]}$ is a bounded operator in $L_\bullet^2(M,E)$. Hence by \cite[Section 4, p. 155]{Ka} the decomposition
\begin{equation}\label{E:dec}
L^2_\bullet(M,E) \, = \, \operatorname{Image} \Pi_{\Delta_{\operatorname{rel}}^\sharp,[0,\lambda]} \oplus \operatorname{Image} (\operatorname{Id} - \Pi_{\Delta_{\operatorname{rel}}^\sharp,[0,\lambda]})
\end{equation}
is a direct sum decomposition into closed subspaces of the Hilbert space $L^2_\bullet(M,E)$. Note that in general the decomposition is not orthogonal with respect to the fixed $L^2$-Hilbert structure unless the Hermitian metric $h^E$ is flat.  

Using the analytic Fredholm theorem, one finds that $\operatorname{Image}\Pi_{\Delta_{\operatorname{rel}}^\sharp,[0,\lambda]} \subset \mathcal{D}(\Delta_{\operatorname{rel}}^\sharp)$ and the projection $\Pi_{\Delta_{\operatorname{rel}}^\sharp,[0,\lambda]}$ commutes with $\Delta_{\operatorname{rel}}^\sharp$. Since $\mathcal{D}(\Delta_{\operatorname{rel}}^\sharp)\subset \mathcal{D}_{\operatorname{min}}$, then, by \eqref{E:dec}, we decompose $\mathcal{D}_{\operatorname{min}}$ into
\[
\mathcal{D}_{\operatorname{min}} \, = \, \mathcal{D}_{\operatorname{min}} \cap \operatorname{Image} \Pi_{\Delta_{\operatorname{rel}}^\sharp,[0,\lambda]} \oplus \mathcal{D}_{\operatorname{min}} \cap \operatorname{Image} (\operatorname{Id} - \Pi_{\Delta_{\operatorname{rel}}^\sharp,[0,\lambda]}).
\]
Denote by $${\mathcal{D}_{\operatorname{min}}}_{,[0,\lambda]} = \mathcal{D}_{\operatorname{min}} \cap \operatorname{Image} \Pi_{\Delta_{\operatorname{rel}}^\sharp,[0,\lambda]}$$ and $${\mathcal{D}_{\operatorname{min}}}_{,(\lambda,\infty)} = \mathcal{D}_{\operatorname{min}} \cap \operatorname{Image} (\operatorname{Id} - \Pi_{\Delta_{\operatorname{rel}}^\sharp,[0,\lambda]}).$$ Since $\nabla_{\operatorname{min}}$ (resp. $\nabla_{\operatorname{min}}^\sharp$) commutes with $\Delta_{\operatorname{rel}}^\sharp$ and hence also with $\Pi_{\Delta_{\operatorname{rel}}^\sharp,[0,\lambda]},$ we find that in fact $({\mathcal{D}_{\operatorname{min}}}_{,\mathcal{I}},{\nabla_{\operatorname{min}}}_{,\mathcal{I}})$ (resp. $({\mathcal{D}_{\operatorname{min}}}_{,\mathcal{I}},{\nabla_{\operatorname{min}}^\sharp}_{,\mathcal{I}})$), $\mathcal{I}=[0,\lambda] \ \text{or} \ (\lambda,\infty)$, where ${\nabla_{\operatorname{min}}}_{,\mathcal{I}}:= \nabla_{\operatorname{min}} |_{{\mathcal{D}_{\operatorname{min}}}_{,\mathcal{I}}}$ (resp. ${\nabla_{\operatorname{min}}^\sharp}_{,\mathcal{I}}:= \nabla_{\operatorname{min}}^\sharp |_{{\mathcal{D}_{\operatorname{min}}}_{,\mathcal{I}}}$), is a subcomplex of the complex $(\mathcal{D}_{\operatorname{min}},\nabla_{\operatorname{min}})$ (resp. $(\mathcal{D}_{\operatorname{min}},\nabla_{\operatorname{min}}^\sharp)$). For each $\lambda \ge 0$, we now have 
\[
(\mathcal{D}_{\operatorname{min}},\nabla_{\operatorname{min}})\, = \, ({\mathcal{D}_{\operatorname{min}}}_{,[0,\lambda]},{\nabla_{\operatorname{min}}}_{,[0,\lambda]})\oplus ({\mathcal{D}_{\operatorname{min}}}_{,(\lambda,\infty)},{\nabla_{\operatorname{min}}}_{,(\lambda,\infty)}).
\]
Analogous decomposition holds for the complex $(\mathcal{D}_{\operatorname{min}},\nabla_{\operatorname{min}}^\sharp)$.
We recall the following proposition, \cite[Corollary 3.16]{Ve1}, \cite[Proposition 3.4]{Su}.
\begin{proposition}\label{P:hocoh}
The subcomplexes $(\mathcal{D}_{\operatorname{min/max},(\lambda,\infty)},\nabla_{\operatorname{min/max},(\lambda,\infty)})$ are acyclic for any $\lambda \ge 0$ and
\begin{equation}\label{E:hocoh}
H^\bullet(\mathcal{D}_{\operatorname{min/max},[0,\lambda]},\nabla_{\operatorname{min/max},[0,\lambda]}) \cong H^\bullet(\mathcal{D}_{\operatorname{min/max}},\nabla_{\operatorname{min/max}}) \cong H^\bullet_{\operatorname{rel}/\operatorname{abs}}(M, E).
\end{equation}
\end{proposition}

We consider the homology groups $H_\bullet(\mathcal{D}_{\operatorname{min}},\nabla_{\operatorname{min}}^\sharp)$ and $H_\bullet({\mathcal{D}_{\operatorname{min}}}_{,[0,\lambda]},{\nabla_{\operatorname{min}}^\sharp}_{,[0,\lambda]})$ now.
By~\eqref{E:chirality} and the fact that $\Gamma \Delta^\sharp_{\operatorname{rel}} \Gamma = \Delta^\sharp_{\operatorname{abs}}$, the chirality operator $\Gamma$ establishes a complex linear isomorphism of the graded complex $({\mathcal{D}_{\operatorname{min}}}_{,[0,\lambda]},{\nabla_{\operatorname{min}}^\sharp}_{,[0,\lambda]})$ to the complex $({\mathcal{D}_{\operatorname{max}}}_{,[0,\lambda]},{\nabla_{\operatorname{max}}}_{,[0,\lambda]})$. Therefore, for $j=0,\cdots,n$, the chirality operator $\Gamma$ establishes a complex linear isomorphism between homology groups and cohomology groups
\begin{equation}\label{E:hocodua}  
\begin{array}{l}
 H_j({\mathcal{D}_{\operatorname{min}}},{\nabla_{\operatorname{min}}^\sharp})
 \cong H^{n-j}({\mathcal{D}_{\operatorname{max}}},{\nabla_{\operatorname{max}}}), \\
 H_j({\mathcal{D}_{\operatorname{min}}}_{,[0,\lambda]},{\nabla_{\operatorname{min}}^\sharp}_{,[0,\lambda]})
\cong H^{n-j}({\mathcal{D}_{\operatorname{max}}}_{,[0,\lambda]},{\nabla_{\operatorname{max}}}_{,[0,\lambda]}). 
\end{array}
\end{equation}
Hence, by \eqref{E:hocoh} and \eqref{E:hocodua}, we have
\begin{equation}\label{E:duaisom}
H_\bullet(\mathcal{D}_{\operatorname{min}},\nabla_{\operatorname{min}}^\sharp) \cong H_\bullet({\mathcal{D}_{\operatorname{min}}}_{,[0,\lambda]},{\nabla_{\operatorname{min}}^\sharp}_{,[0,\lambda]}) \cong H^{n-\bullet}_{\operatorname{abs}}(M,E).
\end{equation}
In particular, we have the following isomorphism
\begin{equation}\label{E:detis}
\operatorname{Det} H_\bullet(\mathcal{D}_{\operatorname{min}},\nabla_{\operatorname{min}}^\sharp)  \cong (\operatorname{Det} H^\bullet_{\operatorname{abs}}(M,E))^{(-1)^n}.
\end{equation}
Hence the Cappell-Miller torsion, cf. \eqref{E:tautor}, can be viewed as an element of the complex line
\begin{align}\label{E:tauminondis1}
T_{\operatorname{min},[0,\lambda]} & := \tau\big(\,{\mathcal{D}_{\operatorname{min}}}_{,[0,\lambda]},{\nabla_{\operatorname{min}}}_{,[0,\lambda]},{\nabla_{\operatorname{min}}^\sharp}_{,[0,\lambda]}\, \big)\nonumber \\
 & \quad \in \operatorname{Det} H^\bullet_{\operatorname{rel}}(M,E)\otimes(\operatorname{Det} H^\bullet_{\operatorname{abs}}(M,E))^{(-1)^{n+1}}.
\end{align}
Using the Poincar\'{e} duality, we also have the isomorphism 
\begin{equation}\label{E:poin}
(\operatorname{Det} H^{n-\bullet}_{\operatorname{abs}}(M,E))^{-1} \cong \operatorname{Det} H^{\bullet}_{\operatorname{rel}}(M,E^*). 
\end{equation}
By~\eqref{E:duaisom} and \eqref{E:poin}, we have
\begin{equation}\label{E:isomo}
\begin{array}{l}
\operatorname{Det} H^\bullet(\mathcal{D}_{\operatorname{min}},\nabla_{\operatorname{min}}) \otimes (\operatorname{Det} H_\bullet(\mathcal{D}_{\operatorname{min}},\nabla_{\operatorname{min}}^\sharp))^{-1} \\
\cong \operatorname{Det} H^\bullet_{\operatorname{rel}}(M,E) \otimes (\operatorname{Det} H^{n-\bullet}_{\operatorname{abs}}(M,E))^{-1} \\
\cong \operatorname{Det} H^\bullet_{\operatorname{rel}}(M,E) \otimes \operatorname{Det} H^\bullet_{\operatorname{rel}}(M,E^*).
\end{array}
\end{equation}
Therefore, the Cappell-Miller torsion can be viewed as an element of the complex line
$$
T_{\operatorname{min},[0,\lambda]}   \in \operatorname{Det} H^\bullet_{\operatorname{rel}}(M,E) \otimes \operatorname{Det} H^\bullet_{\operatorname{rel}}(M, E^*).
$$
In particular, if $n$ is odd, then, by \eqref{E:tauminondis1}, we also have
\begin{equation}\label{E:taumin}
T_{\operatorname{min},[0,\lambda]}  \in \operatorname{Det} H^\bullet_{\operatorname{rel}}(M,E)\otimes \operatorname{Det} H^\bullet_{\operatorname{abs}}(M,E).
\end{equation}

Similarly, we have the Cappell-Miller torsion 
$$
T_{\operatorname{max},[0,\lambda]}   \in \operatorname{Det} H^\bullet_{\operatorname{abs}}(M,E) \otimes \operatorname{Det} H^\bullet_{\operatorname{abs}}(M, E^*).
$$

Denote by ${\Delta_{\operatorname{rel}}^\sharp}_{,(\lambda, \infty)}$ the restriction of ${\Delta_{\operatorname{rel}}^\sharp}$ to $\mathcal{D}(\Delta_{\operatorname{rel}}^\sharp)\cap \operatorname{Image}\big(\operatorname{Id}-\Pi_{{\Delta_{\operatorname{rel}}^\sharp}_{,[0,\lambda]}}\big)$ and, for $0 \le k \le n$, denote by ${\Delta_{\operatorname{rel}}^\sharp}_{,(\lambda, \infty),k}$ the restriction of ${\Delta_{\operatorname{rel}}^\sharp}$ to $\mathcal{D}({\Delta_{\operatorname{rel}}^\sharp}_{,(\lambda,\infty)})\cap L^2_k(M,E)$. Since $\Delta_{\operatorname{rel}}^\sharp$ has the same leading symbol with $\Delta_{\operatorname{rel}}$, the following zeta regularized determinant is well defined, cf. \cite[Theorem 3.9]{Ve1},
\[
\operatorname{Det}({\Delta_{\operatorname{rel}}^\sharp}_{,(\lambda, \infty),k})\, = \, \exp \Big(\, -\frac{\partial}{\partial s}\Big|_{s=0} \operatorname{Tr} [({\Delta_{\operatorname{rel}}^\sharp}_{,(\lambda, \infty),k})^{-s}]    \, \Big).
\]
We also define ${\Delta_{\operatorname{abs}}^\sharp}_{,(\lambda, \infty)}, {\Delta_{\operatorname{abs}}^\sharp}_{,(\lambda, \infty),k}$ and $\operatorname{Det}({\Delta_{\operatorname{abs}}^\sharp}_{,(\lambda, \infty),k})$ in similar ways.

\begin{theorem}
The elements 
\begin{equation}\label{E:taurelabs}
\begin{array}{lcl}
T_{\operatorname{rel/abs}}(\nabla) & := & T_{\operatorname{min/max},[0,\lambda]} \cdot \prod_{k=0}^{n}{\big(\operatorname{Det}(\Delta^\sharp_{\operatorname{rel/abs},(\lambda,\infty),k})\big)^{(-1)^{k+1}k}}\\ 
& & \in  \operatorname{Det} H^\bullet_{\operatorname{rel}/\operatorname{abs}}(M,E) \otimes \operatorname{Det} H^\bullet_{\operatorname{rel}/\operatorname{abs}}(M,E^*)
\end{array}
\end{equation}
are independent of the choice of $\lambda \ge 0$.
\end{theorem}
\begin{proof}
For $0 \le \mu < \lambda < \infty$, one easily sees that
\begin{align}\label{E:detdecompo}
\prod_{k=0}^{n}{\big(\operatorname{Det}(\Delta^\sharp_{\operatorname{rel}/\operatorname{abs},(\mu,\infty),k})\big)^{(-1)^{k+1}k}} & = \prod_{k=0}^{n}{\big(\operatorname{Det}(\Delta^\sharp_{\operatorname{rel}/\operatorname{abs},(\mu,\lambda),k})\big)^{(-1)^{k+1}k}} \nonumber \\
 & \quad \times  \prod_{k=0}^{n}{\big(\operatorname{Det}(\Delta^\sharp_{\operatorname{rel}/\operatorname{abs},(\lambda,\infty),k})\big)^{(-1)^{k+1}k}}.
\end{align}
We also have 
\begin{align}\label{E:decompo1}
(\mathcal{D}_{\min/\max,[0,\lambda]},\nabla_{\min/\max,[0,\lambda]}) & = (\mathcal{D}_{\min/\max,[0,\mu]},\nabla_{\min/\max,[0,\mu]}) \nonumber \\
& \quad \oplus (\mathcal{D}_{\min/\max,[\mu,\lambda]},\nabla_{\min/\max,[\mu,\lambda]})
\end{align}
and
\begin{align}\label{E:decompo2}
(\mathcal{D}_{\min/\max,(\mu,\infty)},\nabla_{\min/\max,(\mu,\infty)}) & = (\mathcal{D}_{\min/\max,(\mu,\lambda)},\nabla_{\min/\max,(\mu,\lambda)}) \nonumber
\\ & \quad \oplus (\mathcal{D}_{\min/\max,(\lambda,\infty)},\nabla_{\min/\max,(\lambda,\infty)}).
\end{align}
By applying the stability property, \cite[Claim C, P.161]{CM}, to \eqref{E:decompo1}, we get
\begin{equation}\label{E:stab}
T_{\min/\max,[0,\mu]} \cdot \prod_{k=0}^{n}{\big(\operatorname{Det}(\Delta^\sharp_{\operatorname{rel/abs},(\mu,\lambda],k})\big)^{(-1)^{k+1}k}} \, = \, T_{\min/\max,[0,\lambda]}.
\end{equation}
Hence, by combining \eqref{E:detdecompo} and \eqref{E:stab}, we obtain the result.
\end{proof}


\begin{definition}
The elements $T_{\operatorname{rel/abs}}(\nabla)$ are called the Cappell-Miller analytic torsions on the complex line $\operatorname{Det} H^\bullet_{\operatorname{rel/abs}}(M,E) \otimes \operatorname{Det} H^\bullet_{\operatorname{rel/abs}}(M,E^*)$.
\end{definition}

\begin{remark}
By \eqref{E:chirality}, we can easily check that
\begin{equation}\label{E:ldlrel}
\Gamma  \Delta^\sharp_{\operatorname{rel}}  \Gamma \, = \, \Delta^\sharp_{\operatorname{abs}}, \qquad \Gamma  \Delta_{\operatorname{rel}}  \Gamma \, = \, \Delta'_{\operatorname{abs}}.
\end{equation}
Hence we have 
\begin{equation}\label{E:releab}
\prod_{k=0}^{n}{\big(\operatorname{Det}(\Delta^\sharp_{\operatorname{rel},(\lambda,\infty),k})\big)^{(-1)^{k+1}\cdot k}}\, = \, \prod_{k=0}^{n}{\big(\operatorname{Det}(\Delta^\sharp_{\operatorname{abs},(\lambda,\infty),k})\big)^{(-1)^{k+1}\cdot k}}
\end{equation}
and, by Hodge isomporphisms
$
\operatorname{ker} \Delta_{\operatorname{rel/abs}} \cong H^\bullet_{\operatorname{rel/abs}}(M,E)
$ 
and 
$
\operatorname{ker} \Delta'_{\operatorname{rel/abs}} \cong H^\bullet_{\operatorname{rel/abs}}(M,E^*),
$
we have
$$
\Gamma \, : \, H^k_{\operatorname{rel/abs}}(M,E) \cong H^{n-k}_{\operatorname{abs/rel}}(M,E^*), \quad k=0,\cdots,n.
$$
Therefore, up to an isomorphism between $\operatorname{Det} H^\bullet_{\operatorname{rel}}(M,E) \otimes \operatorname{Det} H^\bullet_{\operatorname{rel}}(M,E^*)$ and 
$\operatorname{Det} H^\bullet_{\operatorname{abs}}(M,E) \otimes \operatorname{Det} H^\bullet_{\operatorname{abs}}(M,E^*)$, the Cappell-Miller analytic torsions $T_{\operatorname{rel}}(\nabla)$ and $T_{\operatorname{abs}}(\nabla)$ are actually the same.
\end{remark}

\begin{remark}
Denote by 
$$
\begin{array}{l}
\mu_{E,E^*}\, : \, \operatorname{Det} H^\bullet_{\operatorname{rel/abs}}(M,E) \otimes \operatorname{Det} H^\bullet_{\operatorname{rel/abs}}(M, E^*) \\ 
\to \operatorname{Det} \big(\, H^\bullet_{\operatorname{rel/abs}}(M,E) \oplus H^\bullet_{\operatorname{rel/abs}}(M, E^*)\,\big)
\cong \operatorname{Det} H^\bullet_{\operatorname{rel/abs}}(M,E \oplus E^*)
\end{array}
$$ 
the canonical fusion isomorphism, cf. \cite[Subsection 2.3]{BK3}. Then $\mu_{E,E^*}(T_{\operatorname{rel/abs}}(\nabla))$ can be viewed as the Cappell-Miller analytic torsions on the determinant lines $\operatorname{Det} H^\bullet_{\operatorname{rel/abs}}(M,E \oplus E^*)$.
\end{remark}

\begin{remark}\label{R:rem1}
Note that in the case that $E$ is acyclic and the Hermitian metric $h^E$ is flat, these torsions are identical with the square of Ray-Singer analytic torsions on manifolds with boundary, cf. for example \cite{Vi}. In the case that the Hermtian metric $h^E$ is flat, in this case $E^*=E$, these torsions can be viewed as the Ray-Singer torsions on the determinant lines $(\operatorname{Det} H^\bullet_{\operatorname{rel/abs}}(M,E))^2 \cong \operatorname{Det} H^\bullet_{\operatorname{rel/abs}}(M,E \oplus E)$.   
\end{remark}

\begin{theorem}
Let $M$ be an odd dimensional oriented compact Riemannian manifold and $(E,\nabla,h^E)$ be a flat complex vector bundle over $M$, then the Cappell-Miller analytic torsions 
\begin{equation}\label{E:taurelab}
T_{\operatorname{rel/abs}}(\nabla) \, := \, T_{\operatorname{min/max},[0,\lambda]} \cdot \prod_{k=0}^{n}{\big(\operatorname{Det}(\Delta^\sharp_{\operatorname{rel/abs},(\lambda,\infty),k})\big)^{(-1)^{k+1}k}} 
\end{equation}
are independent of the choice of $g^M$ in the interior of $M$.
\end{theorem}
\begin{proof}
Consider a smooth family $g^M(t), t \in \mathbb{R}$ of Riemannian metrics, varying only in the interior of $M$. Since the derivation of \cite[(8.10)]{CM} is of local nature, the variation formula of the Ray-Singer term \cite[(8.10)]{CM} is also valid in our case. By combining this with the variation formula of the algebraic formula in \cite[Lemma 7.1]{CM}, which vanishes near the boundary, we obtain the result. 
\end{proof}
\section{Comparison with the refined analytic torsion}
In this section we compare the Cappell-Miller analytic torsion with the refined analytic torsion on odd dimensional manifolds with boundary. 

Denote by
$$
\big(\, \widetilde{\mathcal{D}},\widetilde{\nabla} \, \big) \, = \, (\mathcal{D}_{\operatorname{min}},\nabla_{\operatorname{min}}) \oplus (\mathcal{D}_{\operatorname{max}},\nabla_{\operatorname{max}})
$$
and 
$$
\widetilde{\Gamma}:= \left( \begin{array}{clcr}  0 & \Gamma
\\ \Gamma & 0 \end{array} \right): L^2_\bullet(M,E) \oplus L^2_\bullet(M,E) \to L^2_\bullet(M,E) \oplus L^2_\bullet(M,E).  
$$
Let
$$
\mathcal{B}=\widetilde{\Gamma}\widetilde{\nabla}+\widetilde{\nabla}\widetilde{\Gamma}, \quad \mathcal{D}(\mathcal{B})=\mathcal{D}(\widetilde{\nabla}) \cap \mathcal{D}(\widetilde{\nabla}^*).
$$
Simple computation shows that 
\begin{equation}\label{E:bdrel}
\mathcal{B}^2 \, = \,  \left( \begin{array}{clcr}  \Delta^\sharp_{\operatorname{rel}} & 0
\\ 0 & \Delta^\sharp_{\operatorname{abs}} \end{array} \right).
\end{equation}
Let $\lambda \ge 0$ be any nonnegative real number. Denote by $\Pi_{\mathcal{B}^2,[0,\lambda]}$ the spectral projection of $\mathcal{B}^2$ onto the eigenspaces of absolute value in $[0,\lambda]$, cf. \eqref{E:specpro}. Hence by \cite[Section 4, p. 155]{Ka} the decomposition
\begin{equation}\label{E:decE}
L^2_\bullet(M,E \oplus E) \, = \, \operatorname{Image} \Pi_{\mathcal{B}^2,[0,\lambda]} \oplus \operatorname{Image} (\operatorname{Id} - \Pi_{\mathcal{B}^2,[0,\lambda]}),
\end{equation}
is a direct sum decomposition into closed subspaces of the Hilbert space $L^2_\bullet(M,E \oplus E)$. Denote by $\mathcal{B}_{(\lambda,\infty)}$ the restriction of $\mathcal{B}$ to $\mathcal{D}(\mathcal{B}) \cap \operatorname{Image} (\operatorname{Id} - \Pi_{\mathcal{B}^2,[0,\lambda]})$ and $\mathcal{B}_{(\lambda,\infty),k}$ the restriction of $\mathcal{B}_{(\lambda,\infty)}$ to $L^2_k(M,E \oplus E)$. Denote by $ \rho_{[0,\lambda]}$ the refined torsion in the sense of \cite[Section 4]{BK3} (see also \cite[(3.8)]{Ve1}). Then the refined analytic torsion $\rho_{\operatorname{an}}(\nabla)$, cf. \cite[(4-6), (4-7), (4-16)]{Ve1}, can be written as
\begin{equation}\label{E:reanto}
\rho_{\operatorname{an}}(\nabla) \, = \, \rho_{[0,\lambda]} \cdot \prod_{k=0}^{n}{\big(\operatorname{Det}(\mathcal{B}^2_{(\lambda,\infty),k})\big)^{ (-1)^{k+1}\cdot k/2}} \cdot e^{\pi i(\eta(\widetilde{\nabla})- \operatorname{rank}\cdot \eta(\mathcal{B}_{\operatorname{trivial}}))},
\end{equation}
where $\eta(\widetilde{\nabla})- \operatorname{rank}\cdot \eta(\mathcal{B}_{\operatorname{trivial}})$ is the $\rho$-invariant for the operator $\mathcal{B}$, restricted to even forms, cf. \cite{Ve1}.
The refined torsion $\rho_{[0,\lambda]}$ 
is an element of the determinant line
\[
\rho_{[0,\lambda]} \in \operatorname{Det}(H^\bullet_{\operatorname{rel}}(M,E) \oplus H^\bullet_{\operatorname{abs}}(M,E)).
\]

For $j=0,\cdots,n$, the chirality operator $\Gamma$ defines an action on cohomology $H^j_{\operatorname{rel/abs}}(M,E)$, also denote by $\Gamma$, and let the sets $\{ e_j \}$ and $\{ \Gamma e_{n-j} \}$ be the bases for $H^j_{\operatorname{rel}}(M,E)$ and $H^j_{\operatorname{abs}}(M,E)$ respectively. As a consequence of the construction of $T_{\min,[0,\lambda]}$ and the choice of the sign in \eqref{E:tautor}, the Cappell-Miller analytic torsion can be described as
\begin{equation}\label{E:ekba}
T_{\min,[0,\lambda]}=\Big(\bigotimes_{j=0}^{n}{[e_j]^{(-1)^j}}\Big) \otimes \Big(\bigotimes_{j=0}^{n}{[\Gamma e_{n-j}]^{(-1)^j}}\Big) \in \operatorname{Det} H^\bullet_{\operatorname{rel}}(M,E)\otimes \operatorname{Det} H^\bullet_{\operatorname{abs}}(M,E).
\end{equation}

Denote by $\mu_{(M,E)}$ the fusion isomorphism for graded vector spaces \cite[(2.18)]{BK3}
$$
\mu_{(M,E)} : \operatorname{Det} H^\bullet_{\operatorname{rel}}(M,E)\otimes \operatorname{Det} H^\bullet_{\operatorname{abs}}(M,E) \cong \operatorname{Det}(H^\bullet_{\operatorname{rel}}(M,E) \oplus H^\bullet_{\operatorname{abs}}(M,E)).
$$
then, by \cite[P.29]{Ve2} and \eqref{E:ekba}, we have 
\begin{equation}\label{E:tmin1}
\mu^{(-1)}_{(M,E)}(\rho_{[0,\lambda]}) \, = \, (-1)^{\nu(M,E)}\cdot T_{\min,[0,\lambda]}.
\end{equation}
Here the integer $\nu(M,E)$ is given by the sum
\begin{align}\label{E:signcoll}
\nu(M,E) & :=  \mathcal{M}(\Omega^\bullet_{\operatorname{rel},[0,\lambda]}(M,E),\Omega^\bullet_{\operatorname{abs},[0,\lambda]}(M,E)) \nonumber \\
& \quad + \mathcal{R}(\Omega^\bullet_{\operatorname{rel},[0,\lambda]}(M,E) \oplus \Omega^\bullet_{\operatorname{abs},[0,\lambda]}(M,E)),
\end{align}
where $\mathcal{M}(\Omega^\bullet_{\operatorname{rel},[0,\lambda]}(M,E),\Omega^\bullet_{\operatorname{abs},[0,\lambda]}(M,E))$ and $\mathcal{R}(\Omega^\bullet_{\operatorname{rel},[0,\lambda]}(M,E) \oplus \Omega^\bullet_{\operatorname{abs},[0,\lambda]}(M,E))$ are defined in \cite[(2-19)]{BK3} and \cite[(4-2)]{BK3}, respectively.

The following theorem gives the relationship between the Cappell-Miller analytic torsion and the refined analytic torsion.
\begin{theorem}\label{T:comp}
Let $(E,\nabla)$ be a flat complex vector bundle over an odd dimensional oriented compact Riemannian manifold $M$. Then
$$
T_{\operatorname{rel}}(\nabla) \, = \,(-1)^{\nu(M,E)} \cdot \mu^{(-1)}_{(M,E)}(\rho_{\operatorname{an}}(\widetilde{\nabla})) \cdot e^{-\pi i(\eta(\widetilde{\nabla})- \operatorname{rank}\cdot \eta(\mathcal{B}_{\operatorname{trivial}}))}, 
$$
where $\mathcal{B}_{\operatorname{trivial}}$ is the odd signature operator of the trivial line bundle over $M$ and $\nu(M,E)$ is defined in \eqref{E:signcoll}.
\end{theorem}
\begin{proof}
By \eqref{E:ldlrel} and \eqref{E:bdrel}, one easily sees that
\begin{equation}\label{E:derel}
\operatorname{Det}(\mathcal{B}^2_{(\lambda,\infty),k}) \, = \, \operatorname{Det}(\Delta^\sharp_{\operatorname{rel},(\lambda,\infty),k})\cdot \operatorname{Det}(\Delta^\sharp_{\operatorname{abs},(\lambda,\infty),k}) \, = \, (\operatorname{Det}(\Delta^\sharp_{\operatorname{rel},(\lambda,\infty),k}))^2.
\end{equation}
Therefore, we have
\begin{equation}\label{E:detbrel}
\prod_{k=0}^{n}{\big(\operatorname{Det}(\mathcal{B}^2_{(\lambda,\infty),k})\big)^{(-1)^{k+1}\cdot k/2}} \, = \,\prod_{k=0}^{n}{\big(\operatorname{Det}(\Delta^\sharp_{\operatorname{rel},(\lambda,\infty),k})\big)^{(-1)^{k+1}\cdot k}}.
\end{equation}
Then, by combining \eqref{E:reanto}, \eqref{E:tmin1}, \eqref{E:detbrel} with \eqref{E:taurelab}, we obtain the result.
\end{proof}

\begin{remark}
The comparison theorem of refined analytic torsion and Burghelea-Haller torsion on odd dimensional manifolds with boundary has been obtained in \cite[Theorem 4.1]{Su} by G. Su. Hence, using \cite[Theorem 4.1]{Su} and Theorem \ref{T:comp}, we can also compare the Cappell-Miller analytic torsion with the Burghelea-Haller analytic torsion on odd dimensional manifolds with boundary. 
\end{remark}

\section{Gluing formula for the Cappell-Miller analytic torsion}

In this sectioin we first recall the definitons and facts that we need from \cite{Ve2} and then establish the gluing formula for the Cappell-Miller analytic torsion in the case that the Hermitian metric is flat. We establish the gluing formula \eqref{E:magl} below by combining Theorem~\ref{T:comp} with the gluing formula for refined analytic torsion, \cite[Theorem 10.6]{Ve2}. Recall that the most intricate part of the discussion of the gluing formula for refined analytic torsion, \cite[Theorem 10.6]{Ve2}, is the splitting formula, \cite[Proposition 8.1]{Ve2}, for the refined torsion $\rho_{[0,\lambda]}$ in the special case $\lambda=0$, which is done by a careful analysis of long exact sequences in cohomology and Poincar\'{e} duality on manifolds with boundary. The discussion in \cite{Ve2} does not rely on the gluing formula of S. Vishik in \cite{Vi}, where only the case of trivial representations is treated. Throughout this section, we assume that $\lambda=0$.

\subsection{Setup for the gluing formula}
Let $M=M_1 \cup_N M_2$ be an odd dimensional oriented closed Riemannian manifold, where $N$ is an embedded closed hypersurface of codimension one which seperates $M$ into two pieces $M_1$ and $M_2$ such that $M_j,j=1,2$ are compact Riemannian manifolds with boundary $\partial M_j=N$ and orientations induced from $M$. Suppose that $\rho : \pi_{1}(M) \rightarrow
U(n, {\Bbb C})$ is a unitary representation and $E = {\widetilde M}
\times_{\rho} {\Bbb C}^{n}$ is a flat bundle, where ${\widetilde M}$
is a universal covering space of $M$. We choose a flat connection
$\nabla$ and extend it to a covariant differential on $\Omega^{\bullet}(M, E)$. Assume the metric structures $(g^M,h^E)$ to be product near the hypersurface $N$. More precisely, we identify, using the inward geodesic flow, a collar neighborhood $U \subset M$ of the boundary $\partial M$ diffeomorphically with $(-\epsilon,\epsilon) \times \partial M, \epsilon > 0$, where the hypersurface $N$ is identified with $\{ 0 \} \times N$. The metric $g^M$ is product over the collar neighborhood of $N$, if over $U$ it is given under the diffeomorphism $\phi:U \to (-\epsilon,\epsilon)\times \partial M$ by
\[
\phi_*g^M|_U=dx^2 \oplus g^M|_N.
\] 
The diffeomorphism $U \cong (-\epsilon,\epsilon)\times N$ shall be covered by a bundle isomorphism $\tilde{\phi}:E|_U \to (-\epsilon,\epsilon) \times E|_N$. The Hermitian metric $h^E$ is product near the boundary, if it is preserved by the bundle isomorphism, i.e. if for all $x \in (-\epsilon,\epsilon)$
\[
\tilde{\phi}_*h^E|_{\{x\}\times N} = h^E|_N.
\] 

The restrictions of $(E,\nabla)$ to $M_j,j=1,2$ give rise to the twisted de Rham complexes $(\Omega_0^\bullet(M_j,E),\nabla_j)$. We denote their minimal and maximal extensions by
\[
(\mathcal{D}_{j,\operatorname{min/max}},\nabla_{j,\operatorname{min/max}}),
\]
respectively. Recall that (cf. Theorem \ref{T:hilbert}, \cite[Theorem 4.1]{BL} or \cite[Theorem 3.2]{Ve1}) these complexes are Fredholm and their cohomology groups can be computed from smooth subcomplexes as follows. Consider for $j=1,2$ the natural inclusions $\iota_j:N \hookrightarrow M_j$ and put
$$
(\Omega^\bullet_{\operatorname{min}}(M_j,E),\nabla), \, \Omega^\bullet_{\operatorname{min}}(M_j,E) \, = \, \{ \, \omega \in \Omega^\bullet(M_j,E) | \iota^*(\omega)=0 \, \},$$
$$
(\Omega^\bullet_{\operatorname{max}}(M_j,E),\nabla), \, \Omega^\bullet_{\operatorname{max}}(M_j,E) \, = \, \Omega^\bullet(M_j,E).
$$ 
The operators $\nabla_j$ yield exterior derivatives on $\Omega^\bullet_{\min/\max}(M_j,E)$. The complexes $(\Omega^\bullet_{\operatorname{min/max}}(M_j,E),\nabla_j)$ are, by Theorem \ref{T:hilbert}, smooth subcomplexes of the Fredholm complexes $(\mathcal{D}_{j,\operatorname{min/max}},\nabla_{j,\operatorname{min/max}})$ with
\[
H^\bullet_{\operatorname{rel/abs}}(M_j,E) \, := \, H^\bullet(\Omega^\bullet_{\operatorname{min/max}}(M_j,E),\nabla_j)
\, \cong \, H^\bullet(\mathcal{D}_{j,\min/\max},\nabla_{j,\min/\max}).
\]
Finally, we denote by $\Gamma_j$ the restriction of $\Gamma$ and, correspondingly, $\Delta_{\operatorname{rel/abs},j}$ the restriction of $\Delta_{\operatorname{rel/abs}}$ to $M_j,j=1,2$.

\subsection{Spectrum for the Laplacian on the splitting manifold}
Next we consider a complex, that takes the splitting $M=M_1 \cup_N M_2$ into account. Let $\iota_j:N \hookrightarrow M_j,j=1,2$ be the natural inclusions. Put
\[
\Omega^\bullet(M_1 \# M_2,E) \, := \, \{\, (\omega_1,\omega_2) \in \Omega^\bullet(M_1,E) \oplus \Omega^\bullet(M_2,E)| \iota_1^*\omega_1 = \iota_2^*\omega_2  \, \}. 
\]
Denote the restrictions of the flat connection $\nabla$ to $M_j,j=1,2,$ by $\nabla_j$ and extend the restrictions to operators on the complexes $\Omega^\bullet(M_j,E),j=1,2$. We put further 
\[
\nabla_S(\omega_1,\omega_2)\, := \, (\nabla_1\omega_1,\nabla_2\omega_2).
\]
This operation respects the transmission condition of $\Omega^\bullet(M_1 \# M_2,E)$ and further its square is obviously zero. Therefore $\nabla_S$ turns the graded vector space $\Omega^\bullet(M_1 \# M_2,E)$ into a complex, denote by
\begin{equation}\label{E:complexS}
(\Omega^\bullet(M_1 \# M_2,E),\nabla_S).
\end{equation}
The natural $L^2$-structure on $\Omega^\bullet(M_1,E) \oplus \Omega^\bullet(M_2,E)$, induced by the metric $g^M$ and $h^E$, is defined on any $\omega=(\omega_1,\omega_2),\eta=(\eta_1,\eta_2)$ as follows
\[
<\omega,\eta>_{L^2}:=\sum_{j=1}^2<\omega_j,\eta_j>_{M_j}.
\]
In order to analyze the associated Laplace operators, consider first the operator $\nabla_S^*$, the adjoint to $\nabla_S$, in $\Omega^\bullet(M_1,E) \oplus \Omega^\bullet(M_2,E)$ with domain of definition $\mathcal{D}(\nabla_S^*)$ consisting of elements $\omega=(\omega_1,\omega_2) \in \Omega^\bullet(M_1,E) \oplus \Omega^\bullet(M_2,E)$ such that the respective linear functionals on any $\eta=(\eta_1,\eta_2) \in \Omega^\bullet(M_1 \# M_2,E)$
\[
L_\omega(\eta)=<\omega,\nabla_S \eta>_{L^2}
\]
are cotinuous in $\Omega^\bullet(M_1 \# M_2,E)$ with respect to the natural $L^2$-norm of $\eta$. As a consequence of Stokes's formula we find for such elements $\omega \in \mathcal{D}(\nabla_S^*)$ that the following transimission condition has to hold
\[
\star \omega = (\star \omega_1, \star \omega_2) \in \Omega^\bullet(M_1 \# M_2,E),
\]
where $\star$ also denotes the restriction of the usual Hodge star operator on $M$ to $M_j,j=1,2$. The Laplacian $\Delta_S=\nabla_S^*\nabla_S+\nabla_S\nabla_S^*$ of the complex \eqref{E:complexS} acts on the obvious domain of definition
\[
\mathcal{D}(\Delta_S) =  \{ \omega \in \Omega^\bullet(M_1 \# M_2,E) | \omega \in \mathcal{D}(\nabla_S^*), \nabla_S \omega \in \mathcal{D}(\nabla_S^*), \nabla^*_S\omega \in \Omega^\bullet(M_1 \# M_2,E)  \}. 
\] 
The $\operatorname{Dom}(\Delta_S)$ is defined as the completion of $\mathcal{D}(\Delta_S)$ with respect to the graph topology norm. The Laplacian $\Delta_S$ with domain $\operatorname{Dom}(\Delta_S)$ is a self-adjoint operator in the $L^2$-completion of $\Omega^\bullet(M_1,E) \oplus \Omega^\bullet(M_2,E)$. 

We recall the following theorem \cite[Theorem 5.2]{Ve2}, which was essentially established by S. Vishik in \cite[Proposition 1.1]{Vi}
\begin{theorem}\label{T:sisom}
The generalized eigenforms of the Laplacian $\Delta_S$ and the generalized eigenforms of the Laplacian $\Delta$ associated to the twisted de Rham complex $(\Omega^\bullet(M,E),\nabla)$ coincide.
\end{theorem}

The following corollary is an obvious consequence of Theorem \ref{T:sisom} and the Hodge isomorphism.
\begin{corollary}\label{C:hodgeiso}
The Laplacian $\Delta_S$ on $\operatorname{Dom}(\Delta_S)$ is Fredholm operator and 
\[
H^\bullet(M_1 \# M_2,E) \, := \, H^\bullet(\Omega^\bullet(M_1 \# M_2,E),\nabla_S) \,  \cong \, H^\bullet_{\operatorname{dR}}(M,E).
\]
\end{corollary}

\subsection{Canonical isomorphisms associated to long exact sequences}

Consider the short exact sequences of complexes
$$
0 \rightarrow \ (\Omega^\bullet_{\operatorname{min}}(M_1,E),\nabla_1) \ \stackrel{\alpha}{\rightarrow}\ (\Omega^\bullet(M_1 \# M_2,E),\nabla_S) \ \stackrel{\beta}{\rightarrow} \  (\Omega^\bullet_{\operatorname{max}}(M_2,E),\nabla_2) \ \rightarrow 0,
$$
$$
0 \rightarrow \ (\Omega^\bullet_{\operatorname{min}}(M_2,E),\nabla_2) \ \stackrel{\alpha'}{\rightarrow}\ (\Omega^\bullet(M_1 \# M_2,E),\nabla_S) \ \stackrel{\beta'}{\rightarrow} \  (\Omega^\bullet_{\operatorname{max}}(M_1,E),\nabla_1) \ \rightarrow 0,
$$
where $\alpha(\omega)=(\omega,0),\alpha'(\omega)=(0,\omega)$ and $\beta(\omega_1,\omega_2)=\omega_2,\beta'(\omega_1,\omega_2)=\omega_1$. The exactness at the first and the second complex of both sequences is clear by construction. The surjectivity of $\beta$ and $\beta'$ is clear, since $\Omega^\bullet_{\operatorname{min}}(M_j,E),j=1,2$, consist of smooth differential forms over $M_j$ which are in particular smooth at the boundary. These short exact sequences of complexes induce long exact sequences on cohomology:
$$
\mathcal{H}: \cdots H^k_{\operatorname{rel}}(M_1,E) \stackrel{\alpha^*}{\rightarrow} H^k(M_1 \# M_2,E) \stackrel{\beta^*}{\rightarrow} H^k_{\operatorname{abs}}(M_2,E) 
\stackrel{\delta^*}{\rightarrow} H^{k+1}_{\operatorname{rel}}(M_1,E) \cdots
$$
$$
\mathcal{H}': \cdots H^k_{\operatorname{rel}}(M_2,E) \stackrel{\alpha'^*}{\rightarrow} H^k(M_1 \# M_2,E) \stackrel{\beta'^*}{\rightarrow} H^k_{\operatorname{abs}}(M_1,E) 
\stackrel{\delta'^*}{\rightarrow} H^{k+1}_{\operatorname{rel}}(M_2,E) \cdots
$$
The long exact sequences induce isomorphisms on determinant lines (cf. \cite{Ni}) in a canonical way, cf. \cite[Section 7]{Ve2},
$$
\Psi : \operatorname{Det} H_{\operatorname{rel}}^\bullet(M_1,E) \otimes \operatorname{Det} H_{\operatorname{abs}}^\bullet(M_2,E) \to \operatorname{Det} H^\bullet(M_1 \# M_2,E),
$$
$$
\Psi' : \operatorname{Det} H_{\operatorname{abs}}^\bullet(M_1,E) \otimes \operatorname{Det} H_{\operatorname{rel}}^\bullet(M_2,E) \to \operatorname{Det} H^\bullet(M_1 \# M_2,E).
$$
The product of the Cappell-Miller analytic torsions is an element 
\begin{align*}
 T_{\operatorname{rel}}(\nabla_1) \otimes T_{\operatorname{abs}}(\nabla_2)   \in & \big(\,\operatorname{Det} H^\bullet_{\operatorname{rel}}(M_1,E)\otimes \operatorname{Det} H^\bullet_{\operatorname{abs}}(M_1,E)\,\big) \\ & \quad \otimes \big(\,\operatorname{Det} H^\bullet_{\operatorname{abs}}(M_2,E)\otimes \operatorname{Det} H^\bullet_{\operatorname{rel}}(M_2,E)\,\big).
\end{align*}

In order to reorder the determinant lines in the gluing formula \eqref{E:magl} appropriately, we introduce the flip isomorphism 
\begin{equation}\label{E:flip}
\Phi ((v_1 \otimes w_1) \otimes (w_2 \otimes v_2))=((v_1 \otimes w_2) \otimes (w_1 \otimes v_2)),
\end{equation}
where, for $i=1,2$, $v_i \in \operatorname{Det} H^\bullet_{\operatorname{rel}}(M_i,E)$ and $w_i \in \operatorname{Det} H^\bullet_{\operatorname{abs}}(M_i,E)$.
By combining the gluing formula for refined analytic torsion \cite[Theorem 10.6]{Ve2} in the case $\lambda=0$ and Theorem \ref{T:comp}, we obtain the gluing formula for the Cappell-Miller analytic torsion in the case that $\lambda=0$. We can also prove the following theorem by combining \eqref{E:taurelab}, \eqref{E:ekba}, \eqref{E:flip} with \cite[Corollary 7.5 and Corollary 10.7]{Ve2}.

\begin{theorem}\label{T:gl}
Let $M \, = \, M_1 \cup_N M_2$ be an odd-dimensional oriented closed Riemannian splitting manifold where $M_j, j=1,2$ are compact bounded Riemannian manifolds with boundary $\partial M_j=N$ and orientation induced from $M$. Denote by $(E,\nabla,h^E)$ a complex flat vector bundle induced by a unitary representation $\rho : \pi_1(M) \to U(n,\mathbb{C})$ and $\nabla_i,i=1,2$ the restrictions of the flat connection $\nabla$ to $M_i$. Assume product structure for the metrics and the vector bundle. The canonical isomorphism
\begin{align*}
\Psi \otimes \Psi'  : & \big( \operatorname{Det} H^\bullet_{\operatorname{rel}}(M_1,E) \otimes \operatorname{Det} H^\bullet_{\operatorname{abs}}(M_2,E)\big) 
\\ & \quad  \otimes \big( \operatorname{Det} H^\bullet_{\operatorname{abs}}(M_1,E) \otimes \operatorname{Det} H^\bullet_{\operatorname{rel}}(M_2,E)\big) \\ & \to
 \operatorname{Det} H^\bullet(M,E) \otimes \operatorname{Det} H^\bullet(M,E)
\end{align*}
is induced by the long exact sequences on cohomology groups:
$$
\mathcal{H}: \cdots H^k_{\operatorname{rel}}(M_1,E) \to H^k(M,E)\to H^k_{\operatorname{abs}}(M_2,E) 
\to H^{k+1}_{\operatorname{rel}}(M_1,E) \cdots
$$
$$
\mathcal{H}': \cdots H^k_{\operatorname{rel}}(M_2,E) \to H^k(M,E) \to H^k_{\operatorname{abs}}(M_1,E) 
\to H^{k+1}_{\operatorname{rel}}(M_2,E) \cdots
$$
Then the gluing formula for the Cappell-Miller analytic torsion is given by the following
\begin{equation}\label{E:magl}
\big((\Psi \otimes \Psi') \circ \Phi\big) (T_{\operatorname{rel}}(\nabla_1) \otimes T_{\operatorname{abs}}(\nabla_2)) \, = \, 2^{\chi(N)} T(\nabla),
\end{equation}
where $T(\nabla)$ denotes the Cappell-Miller analytic torsion for $M$ and $\Phi$ is defined above, cf. \eqref{E:flip}.
\end{theorem}

The gluing formula \eqref{E:magl} is established only in the case that the Hermitian metric is flat, i.e. the representation is unitary. Recall that the Cappell-Miller analytic torsion, in the case that the Hermitian metric is flat, can be viewed as the Ray-Singer torsion on the determinant lines $(\operatorname{Det} H^\bullet(M,E))^2 \cong \operatorname{Det} H^\bullet(M,E \oplus E)$, cf. Remark \ref{R:rem1}. It would be interesting to establish the gluing formula in the case that the Hermitian metric is not necessarily flat.

\bibliographystyle{amsplain}

\end{document}